\documentclass[10pt,a4paper]{amsart}
\usepackage{amsmath}
\usepackage{amsfonts}
\usepackage{amssymb}
\usepackage[american]{babel}
\usepackage{stmaryrd}

\title[An addition theorem and maximal zero-sum free sets]{An addition theorem and maximal zero-sum free sets in $\mathbb{Z}/p\mathbb{Z}$.}
\author{\'Eric Balandraud}
\address{U.P.M.C.\\Laboratoire Combinatoire et Optimisation\\175, rue du Chevaleret 75013 Paris}
\email{balandraud@math.jussieu.fr}
\date{}

\newtheorem{defi}{Definition}
\newtheorem{theo}{Theorem}
\newtheorem*{theo*}{Theorem}

\newtheorem{prop}{Proposition}
\newtheorem{rem}{Remark}
\newtheorem{lem}{Lemma}
\newtheorem{cor}{Corollary}

\begin{document}

\begin{abstract} Using the polynomial method in additive number theory, this article establishes a new addition theorem for the set of subsums of a set satisfying $A\cap(-A)=\emptyset$ in $\mathbb{Z}/p\mathbb{Z}$:
\[|\Sigma(A)|\geqslant\min\left\{p,1+\frac{|A|(|A|+1)}{2}\right\}.\]
 The proof is similar in nature to Alon, Nathanson and Ruzsa's proof of the Erd\"os-Heilbronn conjecture (proved initially by Dias da Silva and Hamidoune \cite{DH}). A key point in the proof of this theorem is the evaluation of some binomial determinants that have been studied in the work of Gessel and Viennot. A generalization to the set of subsums of a sequence is derived, leading to a structural result on zero-sum free sequences. As another application, it is established that for any prime number $p$, a maximal zero-sum free set in $\mathbb{Z}/p\mathbb{Z}$ has cardinality the greatest integer $k$ such that
\[\frac{k(k+1)}{2}<p,\]
proving a conjecture of Selfridge from $1976$.
\end{abstract}

\maketitle

\section*{introduction}
Given two subsets $A$ and $B$ of an abelian group, we define their sumset: $A+B=\{a+b|a\in A,\ b\in B\}$, we denote also $a+B$ the sumset $\{a\}+B$. A first important addition theorem was discovered by Cauchy in $1813$ and has been rediscovered a century later by Davenport:

\begin{theo*}{(Cauchy-Davenport \cite{Ca,Da1,Da2})}
Let $p$ be a prime number, $A$ and $B$ be two subsets of $A\subset\mathbb{Z}/p\mathbb{Z}$, then:
\[|A+B|\geqslant\min\left\{p,|A|+|B|-1\right\}.\]
\end{theo*}
This theorem can easily be extended to the sumset of more than two sets: $\left|\sum A_i\right|\geqslant\min\left\{p,\sum (|A_i|-1)+1\right\}$.

Many proofs of the Cauchy-Davenport Theorem have been published and generalizations have been made in abelian groups or in torsion-free groups; Chowla's Theorem \cite{Cho}, Mann's Theorem \cite{Ma}, Kneser's Theorem \cite{Kn1,Kn2,B}, see also \cite{Na}.

Another topic in addition theory consists in investigating the cardinality of the restricted sumset: $A\dot{+}B=\{a+b|a\in A,\ b\in B,\ a\neq b\}$. In $1964$, Erd\"os and Heilbronn made a famous conjecture that became in $1994$ the following theorem by Dias da Silva and Hamidoune:

\begin{theo*}{(Dias da Silva, Hamidoune \cite{DH})} Let $p$ be a prime number and $A\subset \mathbb{Z}/p\mathbb{Z}$. For a natural integer $h$, denote $h^\wedge A=\underbrace{A\dot{+}\dots\dot{+}A}_{h\ \textrm{times}}$ the set of subsums of $h$ pairwise distinct elements of $A$. Then,
\[|h^\wedge A|\geqslant\min\{p,h(|A|-h)+1\}.\]
\end{theo*}

In this article, we focus our interest on the set of all subsums:
\begin{defi} Let $A\subset \mathbb{Z}/p\mathbb{Z}$, we denote its set of subsums by:
\begin{align*}
\Sigma(A)\  & =\left\{\sum_{x\in I}x|\emptyset\subset I\subset A\right\}\\
\intertext{and we also denote its set of non-trivial subsums by:}
\Sigma^*(A) & =\left\{\sum_{x\in I}x|\emptyset\subsetneq I\subset A\right\}.
\end{align*}
A subset $A$ is called a zero-sum free subset if $0\not\in\Sigma^*(A)$.
\end{defi}

The set of subsums of $A$ can also be seen as a sumset $\Sigma(A)=\sum_{\substack{a\in A\\a\neq 0}}\{0,a\}$. Naturally, we have $\Sigma(A)=\Sigma^*(A)\cup\{0\}$ if $A$ is a zero-sum free set and $\Sigma(A)=\Sigma^*(A)$ otherwise. Therefore $|\Sigma(A)|=|\Sigma^*(A)|+1$ if $A$ is a zero-sum free set and $|\Sigma(A)|=|\Sigma^*(A)|$ otherwise.

In $1968$, Olson was the first to prove a lower bound on the cardinality of $\Sigma(A)$ for a set $A$ such that $A\cap(-A)=\emptyset$:

\begin{theo}\label{theoOl}{(Olson \cite{Ol})}
Let $A\subset\mathbb{Z}/p\mathbb{Z}$. Suppose $A\cap(-A)=\emptyset$. Then
\[|\Sigma(A)|\geqslant\min\left\{\frac{p+3}{2},\frac{|A|(|A|+1)}{2}\right\}.\]
\end{theo}
Olson's theorem is slightly more specific; the term $\frac{|A|(|A|+1)}{2}$ can be replaced by $1+\frac{|A|(|A|+1)}{2}$ under some conditions. This quantity seemed to be a natural lower bound.

The main theorem of this paper proves that the lower bound $1+\frac{|A|(|A|+1)}{2}$ essentially holds. We prove the following:

\begin{theo*} Let $p$ be an odd prime number. Let $A\subset\mathbb{Z}/p\mathbb{Z}$, such that $A\cap(-A)=\emptyset$. We have
\begin{align}
|\Sigma(A)|\ & \geqslant\min\left\{p,1+\frac{|A|(|A|+1)}{2}\right\},\\
|\Sigma^*(A)| & \geqslant\min\left\{p,\frac{|A|(|A|+1)}{2}\right\}.
\end{align}
\end{theo*}
Notice that $(1)$ and $(2)$ are independent, neither implies the other.

In the first part of this article, we will focus on some particular binomial determinants, until we give a suitable expression for them. Binomial determinants are minors of Pascal's triangle, their evaluation is closely related to non-intersecting paths in lattices and Young Tableaux. The evaluation of our binomial determinants does only require algebraic lemmas.

The second and main part describes the ideas of the polynomial method and the proof of the main theorem. The principal idea is that a multivariate polynomial of given degree cannot vanish on a too big cartesian product. The proof of the main theorem is similar in nature to the proof of the theorem of Dias-Da-Silva and Hamidoune (Erd\"os-Heilbronn conjecture) that Alon-Nathanson-Rusza gave using the polynomial method, \cite{ANR1,ANR2}, see also \cite{Na}. The binomial determinants from the first part of the article play a key role in this proof.

The last part develops three applications of the main theorem. The first application generalizes the main theorem to sequences (or multisets), and gives a structural result for zero-sum free sequences. As a second application, we define a new constant of a group: the asymmetric critical number that imitates the definition of the critical number. The value of the asymmetric critical number of $\mathbb{Z}/p\mathbb{Z}$ is given. The last application gives a proof to a conjecture made by Selfridge in $1976$ on the cardinality of a maximal zero-sum-free set in $\mathbb{Z}/p\mathbb{Z}$. It is proven that, for any prime number $p$, the cardinality of a maximal zero-sum free set in $\mathbb{Z}/p\mathbb{Z}$ is the greatest integer $k$ such that $\frac{k(k+1)}{2}<p$.

\section{Some Binomial Determinants}

Binomial determinants appear in several mathematical subjects, they are defined as minors of Pascal's triangle. In a seminal article \cite{GV}, Gessel and Viennot developed a combinatorial interpretation of this determinants. Their interpretation is closely related to configurations of non-intersecting paths and Young tableaux. Since then, numerous articles on these determinants have been published.

\begin{defi} Let $0\leqslant a_1<\dots <a_d$ and $0\leqslant b_1<\dots <b_d$ be natural integers, we call their binomial determinant, the determinant

\[\left(
\begin{array}{c}
a_1,\dots,a_d\\
b_1,\dots,b_d\\
\end{array}
\right)=\left|
\begin{array}{cccc}
\binom{a_1}{b_1} & \binom{a_1}{b_2} & \dots & \binom{a_1}{b_d}\\
\binom{a_2}{b_1} & \binom{a_2}{b_2} & \dots & \binom{a_2}{b_d}\\
\vdots           & \vdots           &       & \vdots  \\
\binom{a_d}{b_1} & \binom{a_d}{b_2} & \dots & \binom{a_d}{b_d)}
\end{array}
\right|.\]
\end{defi}

We are interested in the following family of binomial determinants:

For $d$ and $i$ two integers such that $0\leqslant i\leqslant d$, we define the determinant:
\[D_{d,i}=\left(
\begin{array}{ccccccc}
d-1,& d,& \dots,& d-2+i,& d+i,& \dots,& 2d-1\\
0,& 2,& \dots,& 2(i-1),& 2i,& \dots,& 2(d-1)\\
\end{array}
\right).\]

We will rely on the following two lemmas from \cite{GV} and on a generalization of one of them to compute the exact value of $D_{d,i}$.

\begin{lem}{(Gessel, Viennot \cite{GV})}\label{Det1}
If $b_1\neq 0$, then

\[\left(
\begin{array}{c}
a_1,\dots,a_d\\
b_1,\dots,b_d\\
\end{array}
\right)=\frac{a_1\dots a_d}{b_1\dots b_d}\left(
\begin{array}{c}
a_1-1,\dots,a_d-1\\
b_1-1,\dots,b_d-1\\
\end{array}
\right).\]
\end{lem}

\begin{proof}
For any $(i,j)\in[1,d]^2$, we use the egality:
\[\binom{a_i}{b_j}=\frac{a_i}{b_j}\binom{a_i-1}{b_j-1},\]
and factorize each column (and each row) of the determinant by $b_j$ (or by $a_i$).
\end{proof}

\begin{lem}{(Gessel, Viennot \cite{GV})}\label{Det0}
Where the $a_i$ are consecutive integers and $b_1=0$, we have:
\[\left(
\begin{array}{c}
x,x+1,\dots,x+d-1\\
0,b_2,\dots,b_d\\
\end{array}
\right)=\left(
\begin{array}{c}
x,x+1,\dots,x+d-2\\
b_2-1,\dots,b_d-1\\
\end{array}
\right).\]
\end{lem}

\begin{proof}
For any $(i,j)\in[2,d]\times[1,d]$, we substract the column $i-1$ from the column $i$. The first row of the determinant becomes $(1,0,\dots,0)$. Moreover, since we have:
\[\binom{x-1+i}{b_j}-\binom{x-2+i}{b_j}=\binom{x-2+i}{b_j-1},\]
then developing the determinant on the first row gives the result.
\end{proof}

Since in our determinants $D_{d,i}$ the $a_i$ are not consecutive integers, because one of them is missing, we need another lemma generalizing Lemma \ref{Det0}

\begin{lem}\label{Det0bis}
When $b_1=0$, and the $a_i$ are consecutive integer with a missing term, we have:

\begin{align*}
 & \left(
\begin{array}{cccccccc}
x,& x+1,& \dots,& x+i-1,& x+i+1,& x+i+2,& \dots,& x+d-1\\
0,& b_2,& \dots,& b_i,&b_{i+1},& b_{i+2},& \dots,& b_d\\
\end{array}
\right)\\
= & \left(
\begin{array}{ccccccc}
x,& \dots,& x+i-2,& x+i,& x+i+1,& \dots,& x+d-2\\
b_2-1,& \dots,& b_i-1,& b_{i+1}-1,& b_{i+2}-1,& \dots,& b_d-1\\
\end{array}
\right)\\
 & +\left(
\begin{array}{ccccccc}
x,& \dots,& x+i-2,& x+i-1,& x+i+1,& \dots,& x+d-2\\
b_2-1,& \dots,& b_i-1,& b_{i+1}-1,& b_{i+2}-1,& \dots,& b_d-1\\
\end{array}
\right).
\end{align*}
\end{lem}

\begin{proof}
It is the same idea as the proof of Lemma \ref{Det0}, we substract each column from the following column using the equality:
\[\binom{x-1+i'}{b_j}-\binom{x-2+i'}{b_j}=\binom{x-2+i'}{b_j-1},\]
for any $i'\neq i$.

Subtracting the column $i$ from the column $i+1$ will reveal the sum of two binomial columns. Indeed, we have:
\[\binom{x+i+1}{b_j}-\binom{x+i-1}{b_j}=\binom{x+i}{b_j-1}+\binom{x+i-1}{b_j-1}.\]
The determinant splits into two determinants whose first rows are both $(1,0\dots,0)$, and developing both determinants on their first row gives the expected formula.
\end{proof}

First, we give the value of the extreme binomial determinants of the family: $D_{d,0}$ and $D_{d,d}$

\begin{prop}\label{TheDet} For $d\geqslant 1$, we have:
\begin{align*}
D_{d,0}= & \left(
\begin{array}{c}
d,d+1,\dots,2d-1\\
0,2,\dots,2(d-1)\\
\end{array}
\right)=2^{d(d-1)/2},\\
D_{d,d}= & \left(
\begin{array}{c}
d-1,d,\dots,2d-2\\
0,2,\dots,2(d-1)\\
\end{array}
\right)=2^{(d-1)(d-2)/2}.
\end{align*}
\end{prop}

\begin{proof} We can easily compute $D_{1,0}=1$ and $D_{2,0}=\left|\begin{array}{cc} 1 & 1 \\ 1 & 3\end{array}\right|=2$.

Let us consider that $d\geqslant 2$, applying first Lemma \ref{Det0} then Lemma \ref{Det1}, we have:

\begin{align*}
D_{d,0} & =\left(
\begin{array}{c}
d,d+1,\dots,2d-1\\
0,2,\dots,2(d-1)
\end{array}
\right)\\
 & =\left(
\begin{array}{c}
d,d+1,\dots,2d-2\\
1,3,\dots,2d-3
\end{array}
\right)\\
 & =\frac{(2d-2)!}{(d-1)!}\frac{2^{(d-1)}(d-1)!}{(2d-2)!}\left(
\begin{array}{c}
d-1,d,\dots,2d-3\\
0,2,\dots,2d-4
\end{array}\right)\\
 & = 2^{(d-1)}D_{d-1,0}.
\end{align*}

Since $D_{1,0}=1$ and for $d\geqslant 2$, $D_{d,0}=2^{(d-1)}D_{d-1}$, we have $D_{d,0}=2^{d(d-1)/2}$.

We notice that the binomial determinant $D_{d,d}=\left(
\begin{array}{c}
d-1,d,\dots,2d-2\\
0,2,\dots,2(d-1)\\
\end{array}
\right)$ has its last two terms equal to $2d-2$, so the last column of this determinant is $\left(\begin{array}{c} 0\\ \vdots\\ 0\\ 1\end{array}\right)$. Therefore, developing on the last column gives:
\[D_{d,d}=\left(
\begin{array}{c}
d-1,d,\dots,2d-2\\
0,2,\dots,2(d-1)\\
\end{array}
\right)=\left(
\begin{array}{c}
d-1,d,\dots,2d-3\\
0,2,\dots,2(d-2)\\
\end{array}
\right)=D_{d-1,0}=2^{(d-1)(d-2)/2}.\]
\end{proof}

\begin{prop} The binomial determinants $D_{d,i}$ satisfy the recurrence relation: For $d\geqslant 2$ and $i\in[1,d]$, we have:
\[D_{d,i}=2^{d-1}\frac{d-1}{d-2+i}D_{d-1,i-1}+2^{d-1}\frac{d-1}{d-1+i}D_{d-1,i}.\]
\end{prop}

\begin{proof}
It suffices to apply first Lemma \ref{Det0bis} then Lemma \ref{Det1} to $D_{d,i}$:
\begin{align*}
  & \left(
\begin{array}{ccccccc}
d-1,& d,& \dots,& d-2+i,& d+i,& \dots,& 2d-1\\
0,& 2,& \dots,& 2(i-1),& 2i,& \dots,& 2(d-1)\\
\end{array}
\right)\\
= & \left(
\begin{array}{ccccccc}
d-1,& d,& \dots,& d-3+i,& d+i-1,& \dots,& 2d-2\\
1,& 3,& \dots,& 2(i-1)-1,& 2i-1,& \dots,& 2(d-1)-1\\
\end{array}
\right)\\
  & +\left(
\begin{array}{cccccccc}
d-1,& d,& \dots,& d-3+i,& d-2+i,& d+i,& \dots,& 2d-2\\
1,& 3,& \dots,& 2(i-1)-1,& 2i-1,& 2i+1,& \dots,& 2(d-1)-1\\
\end{array}
\right)\\
\intertext{We can now apply Lemma \ref{Det1} to both binomial determinants:}
= &\frac{(2d-2)!}{(d-2)!(d+i-2)}\frac{2^{(d-1)}(d-1)!}{(2d-2)!}\\
  & \times\left(
\begin{array}{ccccccc}
d-2,& d-1,& \dots,& d-4+i,& d+i-2,& \dots,& 2d-3\\
0,& 2,& \dots,& 2(i-1)-2,& 2i-2,& \dots,& 2(d-1)-2\\
\end{array}
\right)\\
 & +\frac{(2d-2)!}{(d-2)!(d+i-1)}\frac{2^{(d-1)}(d-1)!}{(2d-2)!}\\
 & \times\left(
\begin{array}{cccccccc}
d-2,& d-1,& \dots,& d-4+i,& d-3+i,& d+i-1,& \dots,& 2d-3\\
0,& 2,& \dots,& 2(i-1)-2,& 2i-2,& 2i,& \dots,& 2(d-1)-2\\
\end{array}
\right)\\
= & 2^{d-1}\frac{d-1}{d-2+i}D_{d-1,i-1}+2^{d-1}\frac{d-1}{d-1+i}D_{d-1,i}.
\end{align*}
\end{proof}

We will consider the normalized sequence $D'_{d,i}=\frac{D_{d,i}}{2^{d(d-1)/2-i}}$. The previous recurrence relation becomes:
\[D'_{d,i}=2\frac{d-1}{d-2+i}D'_{d-1,i-1}+\frac{d-1}{d-1+i}D'_{d-1,i}.\]

\begin{prop}\label{TheDetprime} For any integers $d\geqslant 1$ and $i\in[0,d]$, we have:
\[D'_{d,i}=\binom{d}{i}+\binom{d-1}{i-1}.\]
\end{prop}

\begin{proof} Proposition \ref{TheDet} proves that this equality holds for $D'_{d,0}=\frac{D_{d,0}}{2^{d(d-1)/2}}=1=\binom{d}{0}+\binom{d-1}{-1}$, and $D'_{d,d}=\frac{D_{d,d}}{2^{d(d-1)/2-d}}=2=\binom{d}{d}+\binom{d-1}{d-1}$.

Suppose that the relation holds for any $(d,i)$ such that $d<d_0$ then for any $i\in[1,d_0-1]$, we have:
\begin{align*}
D'_{d_0,i} & =2\frac{d_0-1}{d_0-2+i}D'_{d_0-1,i-1}+\frac{d_0-1}{d_0-1+i}D'_{d_0-1,i}\\
           & =2\frac{d_0-1}{d_0-2+i}\left(\binom{d_0-1}{i-1}+\binom{d_0-2}{i-2}\right)+\frac{d_0-1}{d_0-1+i}\left(\binom{d_0-1}{i}+\binom{d_0-2}{i-1}\right)\\
           & =2\frac{d_0-1}{d_0-2+i}\left(\frac{(d_0-1)!}{(i-1)!(d_0-i)!}+\frac{(d_0-2)!}{(i-2)!(d_0-i)!}\right)\\
           & \hspace{3cm}+\frac{d_0-1}{d_0-1+i}\left(\frac{(d_0-1)!}{i!(d_0-1-i)!}+\frac{(d_0-2)!}{(i-1)!(d_0-1-i)!}\right)\\
           & = 2\frac{d_0-1}{d_0-2+i}\frac{(d_0-2)!}{(i-2)!(d_0-i)!}\left(\frac{d_0-1}{i-1}+1\right)\\
           & \hspace{3cm}+\frac{d_0-1}{d_0-1+i}\frac{(d_0-2)!}{(i-1)!(d_0-1-i)!}\left(\frac{d_0-1}{i}+1\right)\\
           & =2\frac{(d_0-1)!}{(i-1)!(d_0-i)!}+\frac{(d_0-1)!}{i!(d_0-1-i)!}\\
           & =2\binom{d_0-1}{i-1}+\binom{d_0-1}{i}\\
           & =\binom{d_0}{i}+\binom{d_0-1}{i-1}.\\
\end{align*}
\end{proof}

\begin{rem} Notice that the normalized binomial determinants $D'_{d,i}$ also satisfy the classical recurrence relation of Pascal's triangle:
\[D'_{d,i}=D'_{d-1,i-1}+D'_{d-1,i}.\]
\end{rem}

\begin{theo}\label{TheDetQuiLeFait} Let $d$ and $i$ be integers such that $0\leqslant i\leqslant d$. We have:
\[D_{d,i}=2^{\frac{d(d-1)}{2}-i}\binom{d}{i}\frac{d+i}{d}.\]
\end{theo}

\begin{proof} For $i=0$ or $i=d$, Proposition \ref{TheDet} already gave the expected value. Otherwise $i\in[1,d-1]$ and this is a consequence of Proposition \ref{TheDetprime}. Indeed:
\begin{align*}
D_{d,i} & =2^{\frac{d(d-1)}{2}-i}D'_{d,i}\\
        & =2^{\frac{d(d-1)}{2}-i}\left(\binom{d}{i}+\binom{d-1}{i-1}\right)\\
        & =2^{\frac{d(d-1)}{2}-i}\binom{d}{i}\left(1+\frac{i}{d}\right)\\
        & =2^{\frac{d(d-1)}{2}-i}\binom{d}{i}\frac{d+i}{d}.
\end{align*}
\end{proof}

\section{The discrete Cube theorem}

We will use the polynomial method to prove our main theorem. This method is based on the Combinatorial Nullstellensatz \cite{Al}, or equivalently on the polynomial Lemma. Recently a short proof of the Combinatorial Nullstellensatz has been proposed by Michalek \cite{M}. This method has numerous applications in additive number theory, see for instance \cite{Kar}. As a multivariate polynomial of given degree cannot vanish on a too big cartesian product, it establishes that a set defined on a cartesian product with a polynomial constraint has to be large enough.

Amongst all these applications, Alon, Nathanson and Ruzsa gave a proof of Erd\"os-Heilbronn conjecture in \cite{ANR1,ANR2}, see also \cite{Na}. Their proof uses the properties of some strict Ballot numbers, which are combinatorial quantities related to lattice paths. The proof of our main theorem is similar to their proof where the binomial determinants $D_{d,i}$ take the role of the strict ballot number. Notice that the strict ballot numbers used in the proof from Alon, Nathanson and Ruzsa have an interpretation as binomial determinants as well.

For our purpose, we will rely on the following formulation of the polynomial Lemma:

\begin{theo}{(Alon, Nathanson, Ruzsa \cite{ANR2})}\label{ANR} Let $p$ be a odd prime number. Let $R(X_0,\dots,X_{d-1})$ be a multivariate polynomial over $\mathbb{Z}/p\mathbb{Z}$. Let $A_0$, $A_1$,\dots,$A_{d-1}$ be non-empty subsets of $\mathbb{Z}/p\mathbb{Z}$, with $|A_i|=k_i$ and define $m=\sum_{i=0}^{d-1}(k_i-1)-deg(R)$. If the coefficient of $\prod_{i=0}^{d-1}X_i^{k_i-1}$ in the polynomial
\[\left(X_0+\dots+X_{d-1}\right)^mR(X_0,\dots,X_{d-1})\]
is not zero then
\[\left|\left\{a_0+a_1+\dots+a_{d-1}|a_i\in A_i\ \textrm{and}\ R(a_0,\dots,a_{d-1})\neq 0\right\}\right|\geqslant m+1.\]
\end{theo}

An interesting recent specification of this theorem is:

\begin{theo}{(Liu, Sun \cite{LSUN})}\label{LiuSun}
 Let $k,m,n$ be positive integers with $k>m(n-1)$ and let $\mathbb{F}$ be a field of characteristic $p$ where $p$ is zero or greater than $K=(k-1)n-(m+1)\binom{n}{2}$. Let $A_1,\dots,A_n$ be subsets of $\mathbb{F}$ for which
\[|A_n|=k\ \textrm{and}\ |A_{i+1}|-|A_i|\in\{0,1\}.\]

Let $P_1(X),\dots,P_n(X)\in\mathbb{F}[x]$ be monic and of degree $m$. Then, we have
\[\left|\left\{a_1+\dots+a_n|a_i\in A_i,\ \textrm{and}\ P_i(a_i)\neq P_j(a_j)\ \textrm{if}\ i\neq j\right\}\right|\geqslant K+1.\]
\end{theo}

The following lemma will state the polynomial involved in the proof of our main theorem. To express this polynomial, we extend the definition of binomial determinants to any couple of $d$-uplets of natural integers $(a_1,\dots,a_d)$ and $(b_1,\dots,b_d)$ by the same determinant formula. (The new cases of this generalization would only give determinants that are zero, or equal to, up to a factor $\pm 1$, a previously defined binomial determinant.)

\begin{lem}\label{devlp} In any field of characteristic $p$ ($p=\infty$ possibly), let $t$ and $d$ be integers with $t<p$.

We have:

\begin{align*}
  & \left(X_0+\dots+X_{d-1}\right)^t\left(\prod_{0\leqslant i<j\leqslant d-1}\left(X_j^2-X_i^2\right)\right)\\
= & \sum_{\substack{(b_0,\dots,b_{d-1})\\ \sum_{i=0}^{d-1} b_i=t+d(d-1)\\\max\{b_i\}<p}}\frac{t!\prod_{i=0}^{d-1} (2i)!}{\prod_{i=0}^{d-1} b_i!}\left(
\begin{array}{c}
b_0,b_1,\dots,b_{d-1}\\
0,2,\dots,2(d-1)\\
\end{array}
\right)\prod_{i=0}^{d-1} X_i^{b_i}.
\end{align*}
\end{lem}

\begin{proof}

Let us start with the left hand side of this equality:
\[L=\left(X_0+\dots+X_{d-1}\right)^t\left(\prod_{0\leqslant i<j\leqslant d-1}\left(X_j^2-X_i^2\right)\right).\]
The first factor can be developed using the multinomial Theorem. The second factor is the VanderMonde determinant of the $X_i^2$ and can be developed as well.
\[\prod_{0\leqslant i<j\leqslant d-1}\left(X_j^2-X_i^2\right)=
\begin{array}{|cccc|}
1 & X_0^2 & \dots & X_0^{2(d-1)}\\
1 & X_1^2 & \dots & X_1^{2(d-1)}\\
\vdots & \vdots & & \vdots \\
1 & X_{d-1}^2 & \dots & X_{d-1}^{2(d-1)}
\end{array}=\sum_{\sigma\in\mathfrak{S}_d}sign(\sigma)\prod_{i=0}^{d-1} X_i^{2\sigma(i)},\]
where $\mathfrak{S}_d$ is the set of permutations of $[0,d-1]$ and $sign(\sigma)$ is the signature of $\sigma$.

Therefore, we have:
\begin{align*}
L = & \left(\sum_{\substack{(t_0,\dots,t_{d-1})\\\sum_{i=0}^{d-1} t_i=t}}\frac{t!}{\prod_{i=0}^{d-1} t_i!}\prod_{i=0}^{d-1} X_i^{t_i}\right)\left(\sum_{\sigma\in\mathfrak{S}_d}sign(\sigma)\prod_{i=0}^{d-1} X_i^{2\sigma(i)}\right)\\
= & t!\sum_{\sigma\in\mathfrak{S}_d}sign(\sigma)\sum_{\substack{(t_0,\dots,t_{d-1})\\\sum_{i=0}^{d-1} t_i=t}}\frac{1}{\prod_{i=0}^{d-1} t_i!}\prod_{i=0}^{d-1} X_i^{t_i+2\sigma(i)}.
\end{align*}

For any $(t_0,\dots,t_{d-1})$ such that $\sum_{i=0}^{d-1} t_i=t$ and any $\sigma\in\mathfrak{S}_d$, let us consider $(b_0,\dots,b_{d-1})$ such that $b_i=t_i+2\sigma(i)$, then $\sum_{i=0}^{d-1} b_i=t+d(d-1)$. Moreover, for any $(b_0,\dots,b_{d-1})$ such that $\sum_{i=0}^{d-1} b_i=t+d(d-1)$ and any $\sigma\in\mathfrak{S}_d$, such that forany $i\in [0,d-1]$ $0\leqslant b_i-2\sigma(i)\leqslant t$, there exists an unique $(t_0,\dots,t_{d-1})$ such that $\sum_{i=0}^{d-1} t_i=t$ and $b_i=t_i+2\sigma(i)$. Thus:

\begin{align*}
L = & t!\sum_{\sigma\in\mathfrak{S}_d}sign(\sigma)\sum_{\substack{(b_0,\dots,b_{d-1})\\0\leqslant b_i-2\sigma(i)\leqslant t\\\sum_{i=0}^{d-1} b_i=t+d(d-1)}}\frac{1}{\prod_{i=0}^{d-1} (b_i-2\sigma(i))!}\prod_{i=0}^{d-1} X_i^{b_i}.
\end{align*}

Moreover, let $(b_0,\dots,b_{d-1})$ and $\sigma\in\mathfrak{S}_d$ such that $\sum_{i=0}^{d-1} b_i=t+d(d-1)$, if $b_i<p$, we can write $\frac{1}{(b_i-2\sigma(i))!}=\frac{(2\sigma(i))!}{b_i!}\binom{b_i}{2\sigma(i)}$. If there exists an index $i_0$ such that $b_{i_0}-2\sigma(i_0)<0$, then $\binom{b_i}{2\sigma(i)}=0$. If for any $i\in[0,d-1]$, we have $b_{i}-2\sigma(i)\geqslant 0$, suppose there exists an index $i_1$ such that $b_{i_1}\geqslant p$, then either $b_{i_1}-2\sigma(i_1)>t$, which would imply that $\sum_{i=0}^{d-1}b_i-d(d-1)>t$, which is impossible, or $b_{i_1}-2\sigma(i_1)\leqslant t<p$ and $\binom{b_{i_1}}{2\sigma(i_1)}=0\pmod{p}$. Therefore, we can write:

\begin{align*}
L = & t!\sum_{\substack{(b_0,\dots,b_{d-1})\\\sum_{i=0}^{d-1} b_i=t+d(d-1)\\\max\{b_i\}<p}}\sum_{\sigma\in\mathfrak{S}_d}sign(\sigma)\frac{\prod_{i=0}^{d-1} (2\sigma(i))!}{\prod_{i=0}^{d-1} b_i!}\prod_{i=0}^{d-1} \binom{b_i}{2\sigma(i)}\prod_{i=0}^{d-1} X_i^{b_i}\\
= & \sum_{\substack{(b_0,\dots,b_{d-1})\\\sum_{i=0}^{d-1} b_i=t+d(d-1)\\\max\{b_i\}<p}}\frac{t!\prod_{i=0}^{d-1} (2i)!}{\prod_{i=0}^{d-1} b_i!}\left(\sum_{\sigma\in\mathfrak{S}_d}sign(\sigma)\prod_{i=0}^{d-1} \binom{b_i}{2\sigma(i)}\right)\prod_{i=0}^{d-1} X_i^{b_i}.
\end{align*}

In the coefficient of the monomial $\prod_{i=0}^{d-1} X_i^{b_i}$ we recognize the binomial determinant:

\[\sum_{\sigma\in\mathfrak{S}_d}sign(\sigma)\prod_{i=0}^{d-1} \binom{b_i}{2\sigma(i)}=
\left|
\begin{array}{cccc}
\binom{b_0}{0} & \binom{b_0}{2} & \dots & \binom{b_0}{2(d-1)}\\
\binom{b_1}{0} & \binom{b_1}{2} &       &   \\
\vdots         &                &       &   \\
\binom{b_{d-1}}{0} & \binom{b_{d-1}}{2} & \dots & \binom{b_{d-1}}{2(d-1)}
\end{array}
\right|=\left(
\begin{array}{c}
b_0,b_1,\dots,b_{d-1}\\
0,2,\dots,2(d-1)\\
\end{array}
\right).\]
\end{proof}

To ease the reading of the proof of the main theorem, we isolate in the following lemma the relevant part of the polynomial Lemma that will be needed.

\begin{lem}\label{genesum} In any field $\mathbb{F}$ of characteristic $p$ ($p=\infty$ eventually), let $A_1$, \dots, $A_d$ be subsets of $\mathbb{F}$ of cardinality $|A_i|=k_i$.

If we have:
\begin{itemize}
\item $\sum_{i=1}^d(k_i-1)-d(d-1)<p$,
\item $2d<p$
\item and if the binomial determinant\\
$\left(
\begin{array}{c}
k_1-1,k_2-1,\dots,k_d-1\\
0,2,\dots,2(d-1)\\
\end{array}
\right)$ is not zero modulo $p$,
\end{itemize}
then the set
\[C=\left\{a_1+a_2+\dots+a_d|a_i\in A_i\ \textrm{and}\ a_i\neq\pm a_j\right\}\]
has cardinality
\[|C|\geqslant\sum_{i=1}^d(k_i-1)-d(d-1)+1.\]
\end{lem}

\begin{proof}
Let us first remark that $(a_1,\dots,a_d)$ satisfy the condition $a_i\neq\pm a_j$ if and only if $\prod_{1\leqslant i<j\leqslant d}\left(a_j^2-a_i^2\right)\neq 0$.
We define $t=\sum_{i=1}^d(k_i-1)-d(d-1)$ and we consider the polynomial:
\[P(X_1,\dots,X_d)=\left(\prod_{1\leqslant i<j\leqslant d}\left(X_j^2-X_i^2\right)\right)\left(X_1+\dots+X_d\right)^t.\]
This polynomial has degree $t+d(d-1)=\sum_{i=1}^d(k_i-1)$ and Lemma \ref{devlp} asserts that the coefficient of the monomial $X_1^{k_1-1}\dots X_d^{k_d-1}$ is
\[\frac{t!\prod_{i=1}^d (2i)!}{\prod_{i=1}^d (k_i-1)!}\left(
\begin{array}{c}
k_1-1,k_2-1,\dots,k_d-1\\
0,2,\dots,2(d-1)\\
\end{array}
\right),\]
which is not zero in the field $\mathbb{F}$, since $t<p$, $2d<p$ and $\left(
\begin{array}{c}
k_1-1,k_2-1,\dots,k_d-1\\
0,2,\dots,2(d-1)\\
\end{array}
\right)$ is not zero modulo $p$. Therefore, Theorem \ref{ANR} asserts that $|C|\geqslant t+1$.
\end{proof}

We can now prove the main theorem.

\begin{theo}\label{Cuberules} Let $p$ be an odd prime number. Let $A\subset\mathbb{Z}/p\mathbb{Z}$, such that $A\cap(-A)=\emptyset$. We have
\begin{align}
|\Sigma(A)|\ & \geqslant\min\left\{p,1+\frac{|A|(|A|+1)}{2}\right\},\\
|\Sigma^*(A)| & \geqslant\min\left\{p,\frac{|A|(|A|+1)}{2}\right\}.
\end{align}
\end{theo}

\begin{proof}

Let us consider a set $A\subset\mathbb{Z}/p\mathbb{Z}$, such that $A\cap(-A)=\emptyset$ and let $|A|=d$. Necessarily, we have $d\leqslant\frac{p-1}{2}$.
Since $p$ is odd, $2$ is invertible modulo $p$, we denote $A=\{2a_1,\dots,2a_d\}$. We have:
\[\Sigma(A)=\sum_{i=1}^d\{0,2a_i\}=\left(\sum_{i=1}^da_i\right)+\sum_{i=1}^d\{-a_i,a_i\},\]
with $a_i\neq\pm a_j$. Therefore, $\sum_{i=1}^d\{-a_i,a_i\}$ and $\Sigma(A)$ have same cardinality.

Let us first consider that $\frac{d(d+1)}{2}<p$.

We consider the sets:
\[\begin{array}{ll}
A_1       & =\left\{a_1,\dots, a_d,-a_1\right\}\\
A_2       & =\left\{a_1,\dots, a_d,-a_1,-a_2\right\}\\
\ \vdots  &  \hspace{4cm} \ddots\\
A_{d-1}   & =\left\{a_1,\dots, a_d,-a_1,\hspace{0.6cm}\dots\hspace{0.6cm},-a_{d-1}\right\}\\
A_d       & =\left\{a_1,\dots, a_d,-a_1,\hspace{0.6cm}\dots\hspace{0.6cm},-a_{d-1},-a_d\right\},
\end{array}\]
and the polynomials for any $i\in[1,d-1]$, $P_i(x)=x^2$.
We check that the set $A_i$ and the polynomials $P_i$ satisfy the conditions of Theorem \ref{LiuSun}. Indeed, since $|A_i|=d+i$, for $i\in[1,d-1]$ we have $|A_{i+1}|-|A_i|=1$ and $|A_d|=2d>2(d-1)$. Finally, the quantity $K$ is:
\begin{align*}
K & =(2d-1)d-(2+1)\binom{d}{2}\\
  & =2d^2-d-3\frac{d^2}{2}+3\frac{d}{2}\\
  & =\frac{d^2}{2}+\frac{d}{2}\\
  & =\frac{d(d+1)}{2}\\
  & <p.
\end{align*}
Therefore, Theorem \ref{LiuSun} states that the set $C=\left\{x_1+\dots+x_n|x_i\in A_i,\ \textrm{and}\ x_i\neq\pm x_j\ \textrm{if}\ i\neq j\right\}$ has cardinality:
\[|C|\geqslant \frac{d(d+1)}{2}+1.\]

Since two terms $x_i$ and $x_j$ of a sum in $C$ cannot be equal or opposite, the elements of $C$ have to be sums of the form $\sum_{i\in I}a_i-\sum_{i\notin I}a_i$ for a set $I\subset[1,d]$. Moreover, for any set $\emptyset\subset I\subset[1,k+1]$, there is $(x_1,\dots,x_d)\in A_1\times\dots\times A_d$ such that $\{x_1,\dots,x_d\}=\{a_i|i\in I\}\cup\{-a_i|i\notin I\}$. Therefore, the set $C$ is exactly the sumset $\sum_{i=1}^d\{-a_i,a_i\}$, and we have $|\Sigma(A)|=|C|\geqslant 1+\frac{d(d+1)}{2}$.

From now, we consider the case $\frac{d(d+1)}{2}\geqslant p$. Without loss of generality, we can suppose that $\frac{d(d-1)}{2}<p$.

We define $i_0=\frac{d(d+1)}{2}-p+1$, then $0< i_0\leqslant d$.

We consider the sets:
\[\begin{array}{ll}
A_1       & =\left\{a_1,\dots, a_d\right\}\\
A_2       & =\left\{a_1,\dots, a_d,-a_1\right\}\\
\ \vdots  &  \hspace{3cm} \ddots\\
A_{i_0}   & =\left\{a_1,\dots, a_d,-a_1,\dots,-a_{i_0-1}\right\}\\
A_{i_0+1} & =\left\{a_1,\dots, a_d,-a_1,\dots,-a_{i_0-1},-a_{i_0},-a_{i_0+1}\right\}\\
\ \vdots  &  \hspace{7cm} \ddots\\
A_d       & =\left\{a_1,\dots, a_d,-a_1,\hspace{1.5cm}\dots\hspace{1cm}\dots\hspace{1cm},-a_d\right\}.
\end{array}\]

We have $|A_i|=d-1+i$ if $i\leqslant i_0$ and $|A_i|=d+i$ if $i>i_0$. Thus
\begin{align*}
\sum_{i=1}^{d}(|A_i|-1)-d(d-1) &=d^2+\frac{d(d+1)}{2}-i_0-d-d(d-1)\\
                               &=d^2+p-1-d^2\\
                               &=p-1<p,
\end{align*}
which is an assumption of Lemma \ref{genesum}. To apply Lemma \ref{genesum} we also need to consider the congruence modulo $p$ of the binomial determinant:
\[D_{d,i_0}=\left(
\begin{array}{ccccccc}
d-1,& d,& \dots,& d-2+i_0,& d+i_0,& \dots,& 2d-1\\
0,& 2,& \dots,& 2(i_0-1),& 2i_0,& \dots,& 2(d-1)\\
\end{array}
\right).\]

Since Theorem \ref{TheDetQuiLeFait} proves that this determinant is exactly $D_{d,i_0}=2^{\frac{d(d-1)}{2}-i_0}\binom{d}{i_0}\frac{d+i_0}{d}$, we can affirm that it has no prime divisor greater than $2d$.  But $d\leqslant\frac{p-1}{2}$, so $p\nmid D_{d,i_0}$. From Lemma \ref{genesum}, if we denote
\[C=\left\{x_1+x_2+\dots+x_{d}|x_i\in A_i\ \textrm{and}\ x_i\neq\pm x_j\right\},\]
we have:
\begin{align*}
|C| & \geqslant \sum_{i=1}^d(|A_i|-1)-d(d-1)+1\\
    & = d^2-d+p-1-d(d-1)+1\\
    & = p.
\end{align*}

Since two terms $x_i$ and $x_j$ of a sum in $C$ cannot be equal or opposite, the elements of $C$ have to be sums of the form $\sum_{i\in I}a_i-\sum_{i\notin I}a_i$ for a set $I\subset[1,d]$. Moreover, for any set $\emptyset\subsetneq I\subset[1,k+1]$, there is $(x_1,\dots,x_d)\in A_1\times\dots\times A_d$ such that $\{x_1,\dots,x_d\}=\{a_i|i\in I\}\cup\{-a_i|i\notin I\}$, and there is no $(x_1,\dots,x_d)\in A_1\times\dots\times A_d$ such that $\{x_1,\dots,x_d\}=\{-a_i|i\in[1,d]\}$.

We deduce that the set $C$ is $\Sigma^*(A)$ up to a translation
\begin{align*}
C & =\left\{\sum_{i\in I}a_i-\sum_{i\notin I}a_i|\emptyset\subsetneq I\subset[1,d]\right\}\\
  & =\left(\sum_{i=1}^d-a_i\right)+\left\{\sum_{i\in I}2a_i|\emptyset\subsetneq I\subset[1,d]\right\}\\
  & =\left(\sum_{i=1}^d-a_i\right)+\Sigma^*(A).
\end{align*}
This proves that $|\Sigma^*(A)|=p$. Necessarily, this also implies that $|\Sigma(A)|=p$.
\end{proof}

In the case $\frac{d(d+1)}{2}<p$, Theorem \ref{LiuSun} was not needed. Indeed, this case could have been proved with Lemma \ref{genesum}, using the fact that the determinant $D_{d,0}=2^{d(d-1)/2}\neq 0\pmod{p}$ from Proposition \ref{TheDet} in the same idea as the other case ($\frac{d(d+1)}{2}\geqslant p$). In the case $\frac{d(d+1)}{2}\geqslant p$, it is possible to use Theorem \ref{LiuSun} with the chosen sets $A_i$ if and only if $i_0=d$, since otherwise the quantity $K=\frac{d(d+1)}{2}\geqslant p$ and $|A_{i_0+1}|-|A_{i_0+1}|=2$, which contradicts two assumptions of Theorem \ref{LiuSun}.

\begin{rem} For any set $A\subset\mathbb{Z}/p\mathbb{Z}$, there exists a couple of sets $A_1$ and $A_2$ such that:
\begin{itemize}
\item $A=A_1\cup A_2$ and $A_1\cap A_2=\emptyset$, 
\item $A_1\cap(-A_1)=\emptyset$,
\item $A_2\cap(-A_2)=\emptyset$,
\item $A_2\subset (-A_1)$. 
\end{itemize}
Since we have the egality $\Sigma(A)=\Sigma(A_1)+\Sigma(A_2)$ and $\Sigma^*(A)=\left(\Sigma^*(A_1)+\Sigma(A_2)\right)\cup\left(\Sigma(A_1)+\Sigma^*(A_2)\right)$, we can deduce from this a lower bound on $\Sigma(A)$ using Cauchy-Davenport Theorem. This idea will be fully developed in the next section where a formulation of Theorem \ref{Cuberules} for sequences is derived.
\end{rem}

\section{Applications}

\subsection{Sets of subsums of a sequence}

Many combinatorial problems concern not only subsets in groups but also multisets (finite collections of elements with repetitions allowed). The first of which is the Erd\"os-Ginzburgh-Ziv Theorem. Numerous topics in additive number theory investigate the structure of sequences whose sets of subsums satisfy various properties, see for instance \cite{CRM,GHK}.

\begin{defi} Let $S=(s_1,\dots,s_n)$ be any finite sequence of elements in $\mathbb{Z}/p\mathbb{Z}$ with repetitions allowed. We denote its set of subsums by:
\begin{align*}
\Sigma(S)\  & =\left\{\sum_{i\in I}s_i|\emptyset\subset I\subset [1,n]\right\}\\
\intertext{and we also denote its set of non-trivial subsums by:}
\Sigma^*(S) & =\left\{\sum_{i\in I}s_i|\emptyset\subsetneq I\subset [1,n]\right\}.
\end{align*}
A sequence $S$ is called a zero-sum free sequence if $0\not\in\Sigma^*(S)$.
\end{defi}

When we consider the subsums of a sequence $S$, the order of its elements is not important, which is why it is often more convenient to write $S=\{(s_i,k_i)\}\subset\mathbb{Z}/p\mathbb{Z}\times\mathbb{N}$, where $k_i$ is the multiplicity of $s_i$ in $S$. We will mainly consider the common multiplicity of a couple $(x,-x)$ in $S$, being the sum of the multiplicities of $x$ and $-x$ in $S$.

We can derive Theorem \ref{Cuberules} for sequences in the following way:

\begin{theo}\label{Sequencerules} Let $p$ be an odd prime number. Let $S=(s_1,\dots,s_n)$ be any finite sequence of elements in $(\mathbb{Z}/p\mathbb{Z})^*$. Denote $(l_1,\dots,l_d)$ all the common multiplicities of $S$, ordered such that $l_1\geqslant l_2\geqslant\dots\geqslant l_d$, then we have:

\begin{align}
|\Sigma(S)|\ & \geqslant\min\left\{p,1+\sum_{i=1}^di.l_i\right\},\\
|\Sigma^*(S)| & \geqslant\min\left\{p,\sum_{i=1}^di.l_i\right\}.
\end{align}
\end{theo}

\begin{proof} Let us consider $A=\{a_1,\dots,a_d\}$ a set of representants of the couples $(x,-x)$ in $S$. (For any $i\in[1,n]$, there exists a unique $j\in[1,d]$ and a unique $\epsilon\in\{\pm 1\}$ such that $s_i=\epsilon a_j$.) Naturally, we have $A\cap(-A)=\emptyset$. Without loss of generality, we can consider that $l_i$ is the common multiplicity of the couple $(a_i,-a_i)$ in $S$. We denote $k_i$ the multiplicity of $-a_i$ in $S$. For $i\in[1,d]$, let us denote $A_i=\{a_1,\dots,a_i\}$, then $\Sigma(A_i)=\sum_{j=1}^i\{0,a_i\}$.
\begin{align*}
\Sigma(S) & =\sum_{i=1}^n\{0,s_i\}\\
          & =\left(\sum_{i=1}^d-k_ia_i\right)+\sum_{i=1}^dl_i\{0,a_i\}\\
          & =\left(\sum_{i=1}^d-k_ia_i\right)+\sum_{i=1}^{l_d}\sum_{j=1}^d\{0,a_j\}+\sum_{i=l_d+1}^{l_{d-1}}\sum_{j=1}^{d-1}\{0,a_j\}+\dots+\sum_{i=l_2+1}^{l_1}\{0,a_1\}\\
          & =\left(\sum_{i=1}^d-k_ia_i\right)+\sum_{i=1}^d(l_i-l_{i+1})\Sigma(A_i),
\end{align*}
where $l_{d+1}=0$.

Therefore, we have:
\[|\Sigma(S)|=\left|\sum_{i=1}^d(l_i-l_{i+1})\Sigma(A_i)\right|.\]

From the Cauchy-Davenport theorem, we have:
\begin{align*}
|\Sigma(S)| & \geqslant \min\left\{p,\sum_{i=1}^d(l_i-l_{i+1})\left(|\Sigma(A_i)|-1\right)+1\right\}\\
\intertext{and from Theorem \ref{Cuberules} we have $|\Sigma(A_i)|-1\geqslant\frac{i(i+1)}{2}=\sum_{j=1}^ij$:}
            & \geqslant \min\left\{p,\sum_{i=1}^d(l_i-l_{i+1})\sum_{j=1}^ij+1\right\}\\
            & =\min\left\{p,\sum_{j=1}^d\sum_{i=j}^dj.l_i-\sum_{j=1}^d\sum_{i=j+1}^dj.l_i+1\right\}\\
            & =\min\left\{p,\sum_{i=1}^di.l_i+1\right\}.
\end{align*}

Similarly, to obtain $(6)$, we observe that either there is a couple $(x,-x)$ such that $x$ and $-x$ are both elements of $S$, which implies that $\Sigma(S)=\Sigma^*(S)$ and $(6)$ holds, or all the $k_i$ are zero.

If all the $k_i$ are zero, since any non-trivial subsum of elements of $S$ contains a non-trivial subsum of elements of $A_d$, we have:

\[\Sigma^*(S)=\sum_{i=1}^{d-1}(l_i-l_{i+1})\Sigma(A_i)+(l_d-1)\Sigma(A_d)+\Sigma^*(A_d).\]

Therefore, we have:
\[|\Sigma^*(S)|=\left|\sum_{i=1}^{d-1}(l_i-l_{i+1})\Sigma(A_i)+(l_d-1)\Sigma(A_d)+\Sigma^*(A_d)\right|.\]

From the Cauchy-Davenport theorem, we have:
\begin{align*}
|\Sigma^*(S)| & \geqslant \min\left\{p,\sum_{i=1}^d(l_i-l_{i+1})\left(|\Sigma(A_i)|-1\right)+(l_d-1)\left(|\Sigma(A_d)|-1\right)+\left(|\Sigma^*(A_d)|-1\right)+1\right\}\\
\intertext{and from Theorem \ref{Cuberules} we have $|\Sigma(A_i)|-1\geqslant\frac{i(i+1)}{2}=\sum_{j=1}^ij$ and $|\Sigma^*(A_d)|-1\geqslant\frac{d(d+1)}{2}-1=\sum_{j=1}^dj-1$:}
            & \geqslant \min\left\{p,\sum_{i=1}^d(l_i-l_{i+1})\sum_{j=1}^ij\right\}\\
            & =\min\left\{p,\sum_{j=1}^d\sum_{i=j}^dj.l_i-\sum_{j=1}^d\sum_{i=j+1}^dj.l_i\right\}\\
            & =\min\left\{p,\sum_{i=1}^di.l_i\right\}.
\end{align*}
\end{proof}

Interestingly, in \cite{Lev} Lev showed that for $S=\{(a_i,l_i)\}$ a sequence of natural integers, the inequality $|\Sigma(A)|\geqslant \sum_{i=1}^di.l_i$ holds without assuming $l_1\geqslant l_2\geqslant ...\geqslant l_d$, but $a_1\leqslant a_2\leqslant ...\leqslant a_d$ instead.

As a application, we give a structural result on zero-sum free sequences in $\mathbb{Z}/p\mathbb{Z}$.

\begin{theo} Let $p$ be an odd prime number and $S$ a zero-sum free sequence of $\mathbb{Z}/p\mathbb{Z}$. For any positive integer $k$, there is in $S$ an element of multiplicity at least $\left\lceil\frac{2}{k}|S|-\frac{2}{k(k+1)}(p-1)\right\rceil$.
\end{theo}

\begin{proof} Since $S$ is a zero-sum free sequence of $\mathbb{Z}/p\mathbb{Z}$, it cannot contain two opposite elements, so the common multiplicity of a couple $(x,-x)$ in $S$ would be either the multiplicity of $x$ or the multiplicity of $-x$.

Let us suppose that every couple $(x,-x)$ has common multiplicity strictly less than $\frac{2}{k}|S|-\frac{2}{k(k+1)}(p-1)$. Since $0\notin\Sigma^*(S)$, we have $|\Sigma^*(S)|<p$, thus Theorem \ref{Sequencerules} asserts that:
\begin{align*}
|\Sigma^*(S)| & \geqslant \sum_{i=1}^di.l_i\\
              & = \sum_{j=1}^d\sum_{i=j}^dl_i.\\
\intertext{Using the equality $|S|=\sum_{i=1}^dl_i$, we have:}
              & \geqslant \sum_{j=1}^{k+1}\left(|S|-\sum_{i=1}^{j-1}l_i\right)\\
              & = (k+1)|S|-\sum_{j=1}^k\sum_{i=1}^{j-1}l_i\\
              & = (k+1)|S|-\sum_{i=1}^k(k+1-i).l_i\\
\intertext{and since $(l_k\leqslant\dots\leqslant l_2\leqslant)l_1<\frac{2}{k}|S|-\frac{2}{k(k+1)}(p-1)$:}
              & >(k+1)|S|-\frac{k(k+1)}{2}\left(\frac{2}{k}|S|-\frac{2}{k(k+1)}(p-1)\right)\\
              & =p-1.
\end{align*}
Since $|\Sigma^*(S)|$ is an integer, it implies that $|\Sigma^*(S)|=p$, which is a contradiction.
\end{proof}

\begin{cor} For $|S|=p-1$ and $k=1$, there exists an element with multiplicity at least $\left\lceil\frac{2}{k}|S|-\frac{2}{k(k+1)}(p-1)\right\rceil=p-1$.
\end{cor}

\begin{cor} In particular, if $|S|\geqslant\frac{p+1}{k}$, then there exists some $g\in\mathbb{Z}/p\mathbb{Z}$ with multiplicity at least $\left\lceil\frac{2}{k(k+1)}\frac{p+2k+1}{k}\right\rceil$ in $S$.
\end{cor}
The case $k=2$ of this last corollary coincides with the prime case of a theorem of Geroldinger and Hamidoune \cite{GH}.

\subsection{An asymmetric critical number}

In additive number theory, an important topic is the determination of the critical number of a group. Let $G$ be an abelian group, we recall the definition of its critical number:
\[cr(G)=\min\left\{l\Big|\forall A\subset G\smallsetminus\{0\}\ \textrm{and}\ |A|\geqslant l,\ \Sigma(A)=G\right\}.\]

The critical number has been first introduced by Erd\"os and Heilbronn \cite{EH} for cyclic group of prime order. Numerous contributions from Olson \cite{Ol}, Diderrich \cite{Di}, Mann and Wou \cite{MW}, Dias da Silva and Hamidoune \cite{DH}, allowed to express the expression of the critical number of any abelian group. The last remaining case has been very recently found by Freeze, Gao and Geroldinger \cite{FGG}. In particular, the value of the critical number of a cyclic group of prime order is $cr\left(\mathbb{Z}/p\mathbb{Z}\right)=\left\lfloor 2\sqrt{p-2}\right\rfloor$, \cite{DH}.

Since $1968$, Olson \cite{Ol} considered the problem of the size of the set of subsums with the extra asymmetric condition $A\cap(-A)=\emptyset$ as assumed in Theorem \ref{theoOl}. The asymmetric hypothesis $A\cap(-A)=\emptyset$ implies that the sets $\{0,a\}$ for $a\in A$ are never arithmetical progressions of same difference. This idea leads to the natural question of the determination of an asymmetric critical number of a group:

\begin{defi}
Let $G$ be an abelian group. If it exists, we call its asymmetric critical number the integer:
\[acr(G)=\min\left\{l\Big|\forall A\subset G\smallsetminus\{0\},\ A\cap(-A)=\emptyset\ \textrm{and}\ |A|\geqslant l,\ \Sigma(A)=G\right\}.\]
\end{defi}

There are groups where this constant cannot be defined. Indeed, in any group of exponant $2$ the condition $A\cap(-A)=\emptyset$ cannot be satisfied. Another example is $\mathbb{Z}/3\mathbb{Z}$, where the condition $A\cap(-A)=\emptyset$ implies $|A|=1$ and $\Sigma(A)\neq\mathbb{Z}/3\mathbb{Z}$. It is not defined for $p=5$ either, because a maximal subset $A\subset\mathbb{Z}/5\mathbb{Z}$ such that $A\cap(-A)=\emptyset$, has cardinality $|A|=2$ and its set of subsums $\Sigma(A)$ is such that $|\Sigma(A)|=4<5$.

From Theorem \ref{Cuberules} $(3)$, we can evaluate $acr\left(\mathbb{Z}/p\mathbb{Z}\right)$ for $p>5$.

\begin{prop} Let $p\geqslant 7$ be a prime number, we have:
\[acr\left(\mathbb{Z}/p\mathbb{Z}\right)=\left\lceil-\frac{1}{2}+\sqrt{2p-\frac{7}{4}}\right\rceil.\]
\end{prop}

\begin{proof} Let $k$  be the greatest number such that $1+\frac{k(k+1)}{2}<p$.
Let $S\subset\mathbb{Z}/p\mathbb{Z}$ be the set $[1,k]$. Since $p>5$, we have $k<\frac{p-1}{2}$, so $S\cap(-S)=\emptyset$. Moreover $\Sigma(S)=[0,\frac{k(k+1)}{2}]\neq\mathbb{Z}/p\mathbb{Z}$. Therefore, $acr\left(\mathbb{Z}/p\mathbb{Z}\right)>k$.

Let us now consider any set $S$ of cardinality  $k'>k$, such that $S\cap(-S)=\emptyset$. Theorem \ref{Cuberules} $(3)$ asserts that $|\Sigma(S)|\geqslant\min\left\{p,1+\frac{k'(k'+1)}{2}\right\}=p$. Therefore, we have $acr\left(\mathbb{Z}/p\mathbb{Z}\right)=k+1$.

The asymmetric critical number is then the minimal number $s$ such that $\frac{s(s+1)}{2}\geqslant p-1$. This last inequality is equivalent to $s^2+s-2(p-1)\geqslant 0$. The equation $s^2+s-2(p-1)=0$ has only one positive root; $-\frac{1}{2}+\sqrt{2p-\frac{7}{4}}$ and has a positive leading coefficient.
\end{proof}

\subsection{An application to a conjecture of Selfridge}

The question of determining the maximal size of a zero-sum free subset in $\mathbb{Z}/p\mathbb{Z}$ is first mentionned by Erd\"os and Heilbronn \cite{EH}. Initially, they conjectured the upper bound: $2\sqrt{p}$. They could prove the upper bound $3\sqrt{6p}$. Few years later in $1968$, Olson \cite{Ol} found a proof that gave the upper bound $2\sqrt{p}$. Erd\"os conjectured a more precise upper bound $\sqrt{2p}$ in $1973$. More specifically, Selfridge conjectured in $1976$ (see for instance \cite{EG,Guy}) that a maximal zero-sum free subset of $\mathbb{Z}/p\mathbb{Z}$ has cardinality $k$, where $k$ is the greatest integer such that:
\[\frac{k(k+1)}{2}<p.\]
A noticable progress has been made in $1996$ by Hamidoune and Z\'emor \cite{HZ}, who proved the upper bound $\sqrt{2p}+5\ln(p)$.

Recently an asymptotic proof of Selfridge's conjecture has been found independently by Deshouillers and Prakash \cite{Des} and by Nguyen, Szemer\'edi and Vu \cite{NSV}:
\begin{theo} Let $p$ be a sufficiently large prime number and $\mathcal{A}$ be a zero-sum free subset of $\mathbb{Z}/p\mathbb{Z}$ with maximal cardinality. Then
\[Card(\mathcal{A})\ \textrm{is the largest integer such that}\ \frac{k(k+1)}{2}\leqslant p+1.\]
\end{theo}

It does not give any contradiction to Selfridge's conjecture since the case $p+1=k(k+1)/2$ (gives $p=(k+2)(k-1)/2$,) holds only for $k=2$, $p=2$ and $k=3$, $p=5$ and the case $p=k(k+1)/2$ holds only for $k=2$, $p=3$. Both articles \cite{Des,NSV} investigate not only the size of a maximal zero-sum free subset, but also the structure of such a set. Precise descriptions of these sets are proved.

We will now prove Selfridge's conjecture for any prime number:

\begin{theo} Let $p$ be a prime number and $A$ be a zero-sum free subset of $\mathbb{Z}/p\mathbb{Z}$ of maximal cardinality, then $|A|$ is the greatest integer $k$ such that:
\[\frac{k(k+1)}{2}<p.\]
\end{theo}

\begin{proof}
For $p=2$, the result is obvious since $\frac{1(1+1)}{2}=1<2<3=\frac{2(2+1)}{2}$. Let $p$ be an odd prime number and $k$ be the greatest integer such that $\frac{k(k+1)}{2}<p$.

For $A=[1,k]$ in $\mathbb{Z}/p\mathbb{Z}$, we have $\Sigma^*(A)=\left[1,\frac{k(k+1)}{2}\right]$ and $0\notin\Sigma^*(A)$. Therefore, a maximal zero-sum free subset has cardinality greater than or equal to $k$.

Suppose that there exists a zero-sum free subset $A$ of cardinality $k'>k$, then Theorem \ref{Cuberules} $(4)$ states:
\begin{align*}
\Sigma^*(A) & \geqslant\min\left(p,\frac{k'(k'+1)}{2}\right)\\
            & =p.
\end{align*}

This implies that $0\in\Sigma^*(A)$, which is a contradiction with $A$ being zero-sum free.
\end{proof}

\end{document}